\documentclass{icmart}
\contact[mgross@dpmms.cam.ac.uk]{DPMMS, Centre for Mathematical Sciences, Wilberforce Road, Cambridge CV3 0WB, United Kingdom}
\contact[siebert@math.uni-hamburg.de]{FB Mathematik, Universit\"at Hamburg, Bundesstra\ss e~55, 20146 Hamburg, Germany}

\usepackage{amscd}
\usepackage{amssymb}
\usepackage{graphicx}
\usepackage{color}
\usepackage[all]{xy}
\usepackage{mathrsfs}
\usepackage{marvosym}
\usepackage{stmaryrd}
\usepackage{srcltx}
\usepackage{hyperref}

\definecolor{shadecolor}{rgb}{1,0.9,0.7}

\newtheorem{theorem}{Theorem}[section]

\newtheorem{proposition}[theorem]{Proposition}

\theoremstyle{definition}
\newtheorem{definition}[theorem]{Definition}

\newtheorem{examples}[theorem]{Examples}

\theoremstyle{remark}
\newtheorem{remark}[theorem]{Remark}

\numberwithin{equation}{section}
\numberwithin{figure}{section}
\renewcommand{\thefigure}{\thesection.\arabic{figure}}

\newcommand {\lfor} {\llbracket}
\newcommand {\rfor} {\rrbracket}

\newcommand{\NN} {\mathbb{N}}
\newcommand{\ZZ} {\mathbb{Z}}

\newcommand{\RR} {\mathbb{R}}
\newcommand{\CC} {\mathbb{C}}

\newcommand{\PP} {\mathbb{P}}
\renewcommand{\AA} {\mathbb{A}}
\newcommand{\GG} {\mathbb{G}}

\newcommand {\shP}  {\mathcal{P}}

\newcommand {\shX}  {\mathcal{X}}
\newcommand {\shY}  {\mathcal{Y}}
\newcommand {\cY}  {\mathcal{Y}}

\newcommand {\fob}  {\mathfrak{b}}

\newcommand {\fom}  {\mathfrak{m}}

\newcommand {\Aff}  {\operatorname{Aff}}

\newcommand {\Aut}  {\operatorname{Aut}}

\newcommand {\Conv} {\operatorname{Conv}}

\newcommand {\Gm} {\GG_m}
\newcommand {\gp}  {{\operatorname{gp}}}

\newcommand {\Hom}  {\operatorname{Hom}}

\newcommand {\im}  {\operatorname{im}}

\newcommand {\Int}  {\operatorname{Int}}

\newcommand {\kk} {\Bbbk}

\newcommand {\Lift}  {\operatorname{Lift}}

\newcommand {\lra}  {\longrightarrow}

\renewcommand {\max} {{\operatorname{max}}}

\renewcommand{\O}  {\mathcal{O}}

\newcommand {\ol} {\overline}

\renewcommand{\P}  {\mathscr{P}}
\newcommand{\oP}  {\overline{\mathscr{P}}}

\newcommand {\Pic}  {\operatorname{Pic}}

\newcommand {\Proj} {\operatorname{Proj}}

\newcommand {\scrR}  {\mathscr{R}}
\newcommand {\scrS}  {\mathscr{S}}
\newcommand {\scrT}  {\mathscr{T}}

\newcommand {\Spec} {\operatorname{Spec}}
\newcommand {\Spf}  {\operatorname{Spf}}

\def\mydate{\ifcase\month \or January\or February\or March\or
April\or May\or June\or July\or August\or September\or October\or 
November\or December\fi \space\number\day,\space\number\year}

\title{Local mirror symmetry in the tropics}

\author[Mark Gross and Bernd Siebert]{Mark Gross and Bernd Siebert
\thanks{M.G.\ was partially supported by NSF grant 1105871 and
1262531.}
}

\begin{document}
\begin{abstract}
We discuss how the reconstruction theorem of \cite{affinecomplex}
applies to local mirror symmetry \cite{CKYZ}. This theorem
associates to certain combinatorial data a degeneration of (log)
Calabi-Yau varieties. While in this case most of the subtleties of
the construction are absent, an important normalization condition
already introduces  rich geometry. This condition guarantees the
parameters of the construction are canonical coordinates in the
sense of mirror symmetry. The normalization condition is also
related to a count of holomorphic disks and cylinders, as
conjectured in \cite{affinecomplex} and partially
proved in
\cite{CLL},\cite{CLT},\cite{CCLT}. We sketch a possible alternative
proof of these counts via logarithmic Gromov-Witten theory.

There is also a surprisingly simple interpretation via rooted trees
marked by monomials, which points to an underlying rich algebraic
structure both in the relevant period integrals and the counting of
holomorphic disks.
\end{abstract}



\maketitle

\section{Introduction}

In \cite{logmirror1}, \cite{affinecomplex}, we proposed a mirror
construction as follows. We begin with a polarized degenerating flat
family $\shX\rightarrow T=\Spec R$ of $n$-dimensional Calabi-Yau
varieties where $R$ is a complete local ring. We consider only
degenerations of a special sort which we term \emph{toric
degenerations}, see \cite{logmirror1}, Def.\ 4.1. Roughly, these are
degenerations for which the central fibre is a union of toric
varieties glued along toric strata, and such that the map
$\shX\rightarrow T$ is locally given by a monomial near the
zero-dimensional strata of the central fibre $X_0$. Associated to
this degeneration we construct the \emph{dual intersection complex}
$(B,\P,\varphi)$, where

(a)~$B$ is an $n$-dimensional integral affine manifold with
singularities (possibly with boundary). In other words, $B$ is a
topological manifold with an open subset $B_0$ with
$\Delta:=B\setminus B_0$ of codimension $\ge 2$, such that $B_0$ has
an atlas of coordinate charts whose transition maps lie in
$\Aff(\ZZ^n)$, the group of integral affine transformations.

(b)~$\P$ is a decomposition of $B$ into convex lattice polyhedra
(possibly unbounded). The singular locus $\Delta$ is typically the
union of codimension two cells of the first barycentric subdivision
of $\P$ not intersecting the interior of a maximal cell of $\P$ nor
containing a vertex of $\P$. There is a one-to-one inclusion
reversing correspondence between elements of $\P$ and toric strata
of $X_0$. The local structure of $\P$ near a vertex is determined by
the fan defining the corresponding irreducible component. The
maximal cells of $\P$ are determined by the toric structure of the
map $\shX\rightarrow T$ near the corresponding zero-dimensional
strata of $X_0$.

(c)~$\varphi$ is a \emph{multi-valued piecewise affine function}.
This is a collection $\{(U_i,\varphi_i)\}$ of $\RR$-valued 
functions $\varphi_i$ on an  open cover $\{U_i\}$ of $B$, with each
$\varphi_i$ piecewise affine linear with respect to the polyhedral
decomposition $\P$, and  $\varphi_i-\varphi_j$ being affine linear
on $U_i\cap U_j$. We assume the slopes of the $\varphi_i$ on cells
of $\P$ to be integral. In this case, $\varphi$ is determined by the
polarization on $\shX$, with local representatives near vertices
given by a piecewise linear function defined by restricting the
polarization to the corresponding irreducible component.

Given this data, we obtain the mirror to the degeneration
$\shX\rightarrow  T$ by reinterpreting $(B,\P,\varphi)$ as the
\emph{intersection complex} of another polarized toric  degeneration
$\shY\rightarrow \Spec\kk\lfor t\rfor$ (in the projective case). 
This time, there is a one-to-one inclusion preserving correspondence
between cells of $\P$ and toric strata of $X_0$, the central fibre
of this new degeneration. The cells of $\P$ are the Newton polytopes
for the polarization restricted to the various strata of $X_0$, and
$\varphi$ is determined by the local toric structure of the map near
zero-dimensional strata.

The prime difficulty in the program lies in reconstructing
$\shY\rightarrow  \Spec\kk\lfor t\rfor$ from the data
$(B,\P,\varphi)$. The main result of \cite{affinecomplex} gives an
algorithm for constructing a \emph{structure} $\scrS$ of walls which
tell us how to construct the degeneration. 

More recently \cite{GHKS} has considered families constructed using
the technology of \cite{affinecomplex} over higher dimensional base
schemes. This  represents a modification of the above procedure. In
the typical example, instead of choosing a fixed polarization on
$\cY$, one chooses a monoid $P$ of polarizations. Let
$Q=\Hom(P,\NN)$ be the dual monoid. Then this data determines a
multi-valued piecewise linear function $\varphi$ taking values in
$Q^{\gp}_{\RR}:=Q^{\gp}\otimes_{\ZZ}\RR$. If $\fom$ is the maximal
monomial ideal of $\kk[Q]$, and $\widehat{\kk[Q]}$ denotes the
completion of $\kk[Q]$ with respect to this ideal, then the
construction gives a family $\shY\rightarrow \Spec\widehat{\kk[Q]}$.

The history of the problem of associating a geometric object (complex
manifold, non-Archimedean space, toric degeneration...) to an integral
affine manifold with singularities began with work of Fukaya \cite{Fuk}.
Fukaya gave a heuristic suggesting that one should be able to construct
the mirror to a K3 surface using objects that look like structures in
two dimensions (in two dimensions, we can think of a structure as just
consisting of a possibly infinite number of unbounded rays). Fukaya observed
that holomorphic disks with boundary on fibres of an SYZ fibration (\cite{SYZ}) 
gave
similar pictures of structures on the mirror side.
In 2004, Kontsevich and Soibelman in \cite{KS} gave the first construction
of a structure, showing how given a two-dimensional affine sphere with
singularities one could construct a consistent structure and from this
structure a non-Archimedean K3 surface. We combined the picture of
toric degenerations
we had been developing independently of the above-mentioned authors
with some ideas from \cite{KS}, allowing us to construct degenerations from 
structures in all dimensions in \cite{affinecomplex}. 

In the first two sections of this paper, we shall illustrate the program
by carrying it out completely for toric Calabi-Yau manifolds, a case
usually referred to as local mirror symmetry \cite{CKYZ}. This particular
case can be viewed as being complementary to the case that the ideas of
\cite{KS} was able to handle.
In the remaining sections, we shall analyze 
enumerative meaning and a tropical interpretation of this construction.

\emph{Acknowledgements}: We would like to thank all people who
influenced our way of thinking about various aspects of our program.
Special thanks go to Mohammed Abouzaid, Paul Hacking, Sean Keel, Diego Matessi
and Rahul Pandharipande. 


\section{Degenerations of toric Calabi-Yau varieties}
\label{Sect: Construction}

Our running example is the construction of the mirror of what is
called ``local $\PP^2$'', the total space $X$ of the canonical
bundle $K_{\PP^2}$ over $\PP^2$. Since $X$ itself is a toric
variety, its anti-canonical divisor $-K_X$ is linearly equivalent to
the sum of toric divisors. There are four toric divisors, the zero
section $S\subset X$, which is the maximal compact subvariety of
$X$, and the preimages $F_0,F_1,F_2$ of the three coordinate lines
in $\PP^2$ under the bundle projection $X\to \PP^2$. Toric methods
show that $S+F_0+F_1+F_2\sim0$ and hence $X$ is a non-compact
Calabi-Yau threefold. The normal bundle $N_{S|X}=\O_{\PP^2}(-3)$ is
determined by the adjunction formula from the Calabi-Yau condition
and it is the dual of an ample line bundle. Hence, by a result of
Grauert \cite{grauert}, any embedded $\PP^2$ in a Calabi-Yau
threefold has an analytic neighbourhood biholomorphic to an analytic
neighbourhood of $S$ in $X$.

For the general description, fix throughout $M=\ZZ^n$,
$M_{\RR}=M\otimes_{\ZZ}\RR$,
$N=\Hom_{\ZZ}(M,\ZZ)$.\footnote{Since $M,N$ will eventually be
treated as data for the mirror side our conventions in this
section are opposite to the usual ones in toric geometry.}
Let $\sigma\subseteq M_{\RR}$ be a compact lattice polytope, and
assume $0\in\sigma$. Define 
\[
C(\sigma)=\{(rm,r)\,|\, m\in\sigma, r\in \RR_{\ge 0}\}
\subseteq M_{\RR}\oplus\RR.
\]
The cone $C(\sigma)$ viewed as a fan defines an affine  toric
variety $X_{\sigma}$. A polyhedral decomposition $\oP$ of $\sigma$
into standard simplices leads to a fan
$\Sigma=\{C(\tau)\,|\,\tau\in\oP\}$ which is a refinement of
$C(\sigma)$. This yields a toric resolution of singularities
$X_{\Sigma}\rightarrow X_{\sigma}$. Assume also that the fan
$\Sigma$ supports at least one strictly convex piecewise linear
function.

\begin{figure}
\center{\input{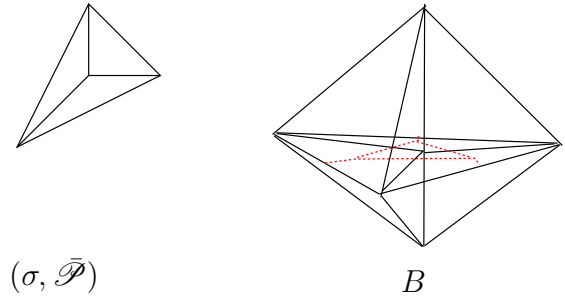}}
\caption{On the left is the initial polytope with its decomposition.
We take the vertices of $\sigma$ to be $(1,0)$, $(0,1)$ and $(-1,-1)$. On
the right is the resulting $B$ with its discriminant locus $\Delta$, indicated
by the dotted line.}
\label{P2localfigure}
\end{figure}

For the case of local $\PP^2$ take $n=2$ and $\sigma=
\Conv\{(1,0),(0,1),(-1,-1)\}$, where $\Conv(S)$ denotes the convex hull
of the set $S$. Then the dual cone $C(\sigma)^\vee$
is $C(\sigma^*)$, the cone over the polar polytope
$\sigma^*$ with vertices $(-1,-1)$, $(2,-1)$, $(-1,2)$. It turns out that
$X_\sigma = \Spec(\CC[C(\sigma)^\vee\cap \ZZ^3])$ is the cyclic
quotient $\AA^2/\ZZ_3$ with $\ZZ_3$ acting diagonally on the
coordinates by multiplication with third roots of unity. Taking the
polyhedral decomposition as shown in Figure~\ref{P2localfigure}
yields for $X_\Sigma$ the blowing up of the origin of $X_\sigma$.
One can show that $X_\Sigma$ is the total space of $K_{\PP^2}$ and
the map to $X_\sigma$ is the contraction of the zero section. Note
also that the projection $C(\sigma)\to M_\RR$ defines a map from
$\Sigma$ to the fan of $\PP^2$, which indeed corresponds to the
bundle projection $X_\Sigma\to \PP^2$.

In general, the map $X_\Sigma\to X_\sigma$ has a reducible
exceptional locus, with one component for each vertex of $\ol\P$
that is not a vertex of $\sigma$, and the explicit description of
the geometry is more complicated.

It turns out that constructing a mirror to $X_{\Sigma}$ does not fit
well with our program. The reason is that $X_\Sigma$ does not seem
to possess a fibration by Lagrangian tori of the kind expected by
mirror symmetry \cite{Examples}. Rather, such a fibration will exist
only after
removal of a hypersurface in $X_{\Sigma}$ that is disjoint from the
exceptional fibre of $X_\Sigma\to X_\sigma$. To run our program we
could give an ad hoc construction of an affine manifold with
singularities derived from the fan $\Sigma$ or write down a toric
degeneration of $X_\Sigma$. The local $\PP^2$ case has been
discussed from the former point of view in \cite{invitation},
Examples~5.1 and 5.2. Since it can be done easily in the present
case we follow the latter method here. This method is motivated by
the construction of toric degenerations of hypersurfaces in toric
varieties in \cite{GBB}.

To exhibit $X_\Sigma$ as an anticanonical hypersurface in a toric
variety we embed the fan $\Sigma$ in $M_\RR\oplus\RR$ as a subfan of
a fan $\tilde\Sigma$ in $M_\RR\oplus\RR^2$. For each maximal cone
$C\in\Sigma$ the fan $\tilde\Sigma$ has two maximal cones\\[-2.5ex]
\[
C_1=
C\times0+\RR_{\ge0} \cdot(0,1,-1),\quad
C_2=
C\times0+\RR_{\ge0} \cdot(0,0,1).
\]
Then $\Sigma$ is the subfan of $\tilde\Sigma$
consisting of cones lying in the hyperplane $M_\RR\oplus\RR\oplus
0\subset M_\RR\oplus\RR^2$. The fan $\tilde\Sigma$ only has two rays
not contained in $\Sigma$, with generators $(0,0,1)$ and $(0,1,-1)$.
The inclusion $M_{\RR}\oplus\RR\oplus 0\subset M_{\RR}\oplus\RR^2$
induces a map of fans from $\Sigma$ to $\tilde\Sigma$, hence an embedding
$j:X_{\Sigma}\hookrightarrow X_{\tilde\Sigma}$ identifying $X_{\Sigma}$
with the closure of the orbit of the subtorus defined by this inclusion
through the distinguished point
(the unit of the toric variety).
Note that the projection to $\RR^2$ maps $\tilde\Sigma$ to the fan
$\Sigma_{\widehat\AA^2}$ of the toric blowing up $\widehat\AA^2$ of
$\AA^2$, with rays generated by $(0,1),(1,0),(1,-1)$. Under this map
the subfan $\Sigma\subset\tilde\Sigma$ maps to the interior ray
$\RR_{\ge0}\cdot(1,0)$. Viewing this interior ray as giving a map
of fans, from the one-dimensional fan defining $\AA^1$ to the two-dimensional
fan defining $\widehat\AA^2$, we obtain an embedding $i:\AA^1\hookrightarrow
\widehat\AA^2$.

We thus obtain a cartesian diagram of
toric morphisms
\[\begin{CD}
X_\Sigma@>j>> X_{\tilde\Sigma}\\
@VpVV@VVqV\\
\AA^1@>i>>  \widehat\AA^2
\end{CD}\\[1ex]
\]
The left vertical arrow is induced by the projection
$M_\RR\oplus\RR\to\RR$, hence is given by the pull-back to
$X_\Sigma$ of the distinguished monomial $x$ on $X_\sigma$
defining the toric boundary as a reduced subscheme. 

Explicitly, write $x,y$ for the toric coordinates on $\AA^2$ and
$\widehat\AA^2=(xu-yv=0)\subset \AA^2\times\PP^1$ for the blowing
up. Then $\im(i)$ is the strict transform of the diagonal $x=y$.
Dehomogenizing $u=1$ or $v=1$ we obtain the usual two coordinate
patches with coordinates $y,v$ and $x,u$ respectively with the
transitions $v=u^{-1}$ and $x=yv$ or $y=xu$. We use the same
notation for the pull-back of $x,y,u,v$ to the corresponding two
types of affine patches with $u\neq 0$ or $v\neq0$ of $X_{\tilde\Sigma}$. 

To describe $X_{\tilde\Sigma}$ let $C\in\Sigma$ be a maximal cone.
Then if $(m,a) \in N\oplus\ZZ$ defines a facet $C'\subset C$, that
is, $(m,a)$ generates an extremal ray of $C^\vee$, the element
$(m,a,a)\in N\oplus\ZZ^2$ defines the facet $C'+\RR_{\ge 0}(0,1,-1)$
of $C_1$. There is only one more facet of $C_1$, namely $C$ itself,
defined by $(0,0,-1)$, and hence 
\[
C_1^\vee= \{(m,a,a)\,|\, (m,a)\in C^{\vee}\}+ \RR_{\ge 0}\cdot(0,0,-1).
\]  
The rays of $C_2^\vee$ are
generated by $(m,a,0)$ for $(m,a)$ an extremal ray of $C^\vee$, and
by $(0,0,1)$, so $C_2^{\vee}=C^{\vee}\times 0+\RR_{\ge 0}(0,0,1)$. 
In either case, we have an identification
\[
\Spec \kk[C_i^\vee\cap (N\oplus\ZZ^2)]=
\Spec \kk[C^\vee\cap N]\times\AA^1
\subset X_\Sigma\times\AA^1.
\]
The toric coordinate for $\AA^1$ is $v= z^{(0,0,-1)}$
for $C_1$ and $u=z^{(0,0,1)}$ for $C_2$. From this description it
is clear that the embedding of $X_\Sigma$ in $X_{\tilde\Sigma}$ is
given by $u=1$ in affine patches with $v\neq0$ and by $v=1$ in the
affine patches with $u\neq 0$.

To write down a degeneration of $X_\Sigma$ to the toric boundary
$\partial X_{\tilde\Sigma} \subset X_{\tilde\Sigma}$ view $u,v$ as
sections of the line bundle $q^*\O(-E)$ where $E\subset \widehat
\AA^2$ is the exceptional curve. Then $X_\Sigma$ is the zero locus
of $s:=u-v$. On the other hand, $xu=yv$ defines a section $s_0$ of
$q^*\O(-E)$ with zero locus $\partial X_{\tilde\Sigma}$. Thus the
hypersurface $\shX\subset X_{\tilde\Sigma}\times\AA^1$ with equation
\[
ts+s_0=0
\]
defines a pencil in $X_{\tilde\Sigma}$ with members $X_\Sigma$ at
$t=\infty$ and with the toric boundary $\partial X_{\tilde\Sigma}$
at $t=0$. Note this pencil is the preimage of the pencil on
$\widehat \AA^2$ defined by the same equations. In particular, by
direct computation $\shX_t$ is completely contained in 
either type of coordinate patch for $t\neq0$. Working in a patch
with $v\neq 0$ we have $s=u-1$, $s_0=xu$ and the equation
\[
0=t s+s_0=t(u-1)+xu= u(t+x)-t
\]
shows $u(t+x)=t\neq0$. Thus $t+x\neq 0$ and $u$ can be eliminated. In
other words, $\shX_t\simeq X_\Sigma\setminus Z_t$ with $Z_t\subset
X_\Sigma$ the hypersurface $x=-t$. Note also that our notation is
consistent in that $x$ indeed descends to the defining equation of
the toric boundary of $X_\sigma$.

It is not difficult to show that $\shX\to\AA^1$ is a toric degeneration.
Indeed, we have already checked that $\shX_0$ is the
toric boundary of $X_{\widetilde{\Sigma}}$. Some harder work shows
that locally near the zero-dimensional strata of $X_0$, the
projection $\shX\rightarrow \AA^1$ is toric. We omit the details,
but this can be done similarly to arguments given in \cite{GBB}.

The dual intersection complex is then easily described along the lines
given in \cite{GBB}, where, for a Calabi-Yau hypersurface in a toric
variety, $B$ was described as the boundary of a reflexive
polytope, with the cones over the faces of the polytope yielding the fan
defining the ambient toric variety.
Topologically, we can write $B\subseteq M_{\RR}\oplus\RR^2$ as
\[
B=\tilde\sigma_1\cup \tilde\sigma_2
\]
where
\begin{align*}
\tilde\sigma_1 = {} & \Conv\big((0,1,-1)\cup (\sigma \times \{(1,0)\})\big),\\
\tilde\sigma_2 = {} & \Conv\big((0,0,1)\cup(\sigma \times \{(1,0)\})\big).
\end{align*}
Note that the support of the fan $\tilde\Sigma$ above is the cone over
$\tilde\sigma_1\cup \tilde\sigma_2$. We then take $\P=\{C\cap B\,|\,C\in
\tilde\Sigma\}$.

Finally, the affine structure on $B$ is defined as follows.
Identify  $\sigma$ with $\sigma\times \{(1,0)\}\subseteq B$, and
take the discriminant locus $\Delta$ to be the union of cells of the
first barycentric subdivision of $\oP$ not containing vertices of
$\oP$, see Figure~\ref{P2localfigure}.  We then define affine charts
as follows. First, we define affine charts $\iota_i: \tilde\sigma_i\setminus
\sigma \hookrightarrow \AA_i$ as the inclusions, 
where $\AA_i$ denotes the affine hyperplane in $M_{\RR}\oplus\RR^2$
spanned by $\tilde\sigma_i$. Second, for each vertex $v\in\oP$, choose a
neighbourhood $U_v$ of $(v,1,0)\in B$. These neighbourhoods can be
chosen so that $U_v\cap U_{v'}=\emptyset$ if $v\not=v'$ and  the two
sets $\tilde\sigma_i\setminus\sigma$ along with the open sets $U_v$ cover
$B\setminus\Delta$. Define a chart $\iota_v:U_v\rightarrow
(M_{\RR}\oplus\RR^2)/\RR(v,1,0)$ via the inclusion followed by the
projection. It is easy to check that these charts give an integral
affine structure. This again precisely follows the procedure for
Calabi-Yau hypersurfaces in toric varieties considered in \cite{GBB}.
This gives rise to the pair $(B,\P)$.

In general, a pair $(B,\P)$ can be described by specifying the
lattice polytopes in $\P$ and specifying a \emph{fan structure} at
each vertex $v$, that is, the identification of a neighbourhood of each
vertex with the  neighbourhood of $0$ in a fan $\Sigma_v$. This
identification gives a one-to-one inclusion preserving 
correspondence between cells of $\P$ containing $v$ and  cones of
$\Sigma_v$, along with integral affine identifications of the
tangent wedges of each cell $\tau\in\P$ containing $v$ with the 
corresponding cone of $\Sigma_v$. These identifications patch
together to give an affine chart in a neighbourhood of the vertex
$v$.

In our example, it is worth describing the fan structure at a vertex
$v\in\sigma$. Since the fan structure at a vertex must be the fan yielding
the corresponding irreducible component of $\shX_0$, 
toric geometry tells us this
fan structure must be given as the quotient fan obtained from
$\tilde\Sigma$ by dividing
out by the ray generated by this vertex. Explicitly,
we use the chart 
\begin{equation}
\label{psivid}
\iota_v:U_v\rightarrow (M_{\RR}\oplus\RR^2)/\RR\cdot(v,1,0)
\cong M_{\RR}\oplus\RR,
\end{equation}
the latter isomorphism given by $(m,r_1,r_2)\mapsto (m-r_1 v,r_2)$.
The fan $\Sigma_v$ can then be described as the fan of tangent wedges
to images of cells $\tau\in\P$ containing $v$.
The set of maximal cones of this fan, described 
as subsets of $M_{\RR}\oplus\RR$, is
\begin{equation}
\label{fanstructureatv}
\{T_v\tau + \RR_{\ge 0}(-v,-1)\,|\,v\in\tau\in\oP_{\max}\}
\cup\{T_v\tau+\RR_{\ge 0}(0,1)\,|\, v\in\tau\in\oP_{\max}\},
\end{equation}
where $T_v\tau$ denotes the tangent wedge to $v\in\tau$ in $M_{\RR}\oplus 0$.
Figure \ref{localA1example} shows some of the fan structures when
$\sigma$ is an interval $[-1,1]$ of length two.

\begin{figure}
\center{\input{localA1example2.pstex_t}}
\caption{}
\label{localA1example}
\end{figure}

One can understand the nature of the singularities of $B$
by studying the local system $\Lambda$ 
of integral vector fields on $B_0$. 
Given integral affine coordinates $y_1,\ldots,y_n$,
$\Lambda$ is locally the family of lattices in the tangent bundle
of $B_0$ generated by $\partial/\partial y_1,\ldots,\partial/\partial y_n$.
If $v,v'
\in\oP$ are adjacent vertices, consider a path $\gamma$ passing from $v$
through $\tilde \sigma_1$ to $v'$ and then through $\tilde\sigma_2$ 
back to $v$. To identify
$\Lambda_v$, we can use the chart \eqref{psivid},
which gives an identification of
$\Lambda_v$ with $M\oplus\ZZ$. It is an
easy exercise to check that parallel transport in $\Lambda$ around $\gamma$
yields a \emph{monodromy transformation} 
\begin{align}
\label{monodromyformula}
\begin{split}
T_{vv'}:\Lambda_{v}&\rightarrow\Lambda_v\\
(m,r)&\mapsto (m+r(v-v'),r)
\end{split}
\end{align}

The final piece of data for the dual intersection complex
$(B,\P,\varphi)$ of $\shX\to\AA^1$ is a multi-valued piecewise
linear function $\varphi$ describing aspects of the K\"ahler
geometry of the situation. 

In the next section, we will build the mirror family over some base
scheme. The natural choice for this base is related to the
K\"ahler cone of $X_\Sigma$, or the Picard group. By toric geometry,
$\Pic(X_\Sigma)$ equals piecewise linear functions on $\Sigma$
modulo linear functions. It will thus be convenient to normalize the
piecewise linear functions as follows. Choose a maximal cell
$\tau\in\oP_{\max}$ which has $0$ as a vertex. Let $P$ be the monoid
of integral convex piecewise linear functions on the fan $\Sigma$
which take the value $0$ on the cone $C(\tau)$. Note that
$P^{\gp}\cong \Pic X_{\Sigma}$. Setting $Q:=\Hom(P,\NN)$,
there is a universal piecewise linear function
$\psi:|\Sigma|\rightarrow Q^{\gp}_{\RR}$, with 
\[
\psi(x)= (P\ni \varphi \mapsto \varphi(x)).
\]
This function is strictly convex in
the sense of \cite{GHK}, Definition 1.12. In the local $\PP^2$ case
normalized piecewise linear functions are determined by
the value at the one remaining vertex of $\sigma$ not contained in
$\tau$ and hence $Q=\NN$.

The multi-valued piecewise linear function $\varphi$ comes from the 
universal piecewise linear function $\psi$ on $|\Sigma|$ by descent
to a quotient fan, or rather from a choice of extension of this
function to $\widetilde{\Sigma}$. This choice can be made by
choosing an element $q\in Q\setminus \{0\}$. While the choice of
$q$ affects the family of polarizations on $X_{\tilde\Sigma}$, it
does not affect the family after restriction to $\shX_t$ for
$t\not=0$. However, it does affect the polarization on $\shX$, and
hence will play some role in the mirror, seen explicitly in
\eqref{Eqn for shY}. We take $\widetilde\psi$ to be the
$Q^{\gp}_{\RR}$-valued piecewise linear extension of $\psi$ which
takes the value $0$ at $(0,1,-1)$ and the value $q$ at $(0,0,1)$.
One can check that this function is strictly convex in the sense of
\cite{GHK}, Definition 1.12. 

We can then construct $\varphi$ from $\widetilde{\psi}$ as follows.
For each $C\in\tilde\Sigma$, let $\tau=C\cap B$ be the corresponding
cell of $\P$. The function $\widetilde{\psi}$ induces a function 
on the quotient fan of $\widetilde\Sigma$ along $C$ (this quotient fan
determining
the fan structure of $B$ along $\tau$) as follows. 
Let $\widetilde\psi_{\tau} \in \Hom(M\oplus\ZZ^2,
Q)$ be a linear extension of  $\widetilde\psi|_{C}$. Then
$\widetilde \psi-\widetilde\psi_{\tau}$ is a piecewise affine
function on $\tilde\Sigma$ vanishing on $C$, hence descending to the
quotient fan of $\tilde\Sigma$ along $C$. 
We take $(\widetilde
\psi-\widetilde\psi_{\tau})|_B$ as a representative of $\varphi$ on
a small open neighbourhood of $\Int(\tau)$ in $B$; this is clearly the
pull-back of the corresponding function on the quotient fan of 
$\tilde\Sigma$ along $C$ under the projection to $(M_{\RR}\oplus\RR^2)/\RR C$.


\section{The mirror degeneration and slab functions}
\label{Sect: mirror}
Having described $(B,\P,\varphi)$ in our example as the dual
intersection complex of a degeneration of the local Calabi-Yau
$X_\Sigma$, we turn to the construction of the mirror, which shall
be a family $\shY\rightarrow \Spec\widehat{\kk[Q]}$ over a generally
higher-dimensional base. 

This family is constructed by constructing
families $\shY_k\rightarrow \Spec\kk[Q]/\fom^{k+1}$ to each order
$k$, giving rise to a formal scheme $\widehat\shY\rightarrow
\Spf\widehat{\kk[Q]}$. As the case at hand will be projective, the
Grothendieck existence theorem gives rise to the desired family.
Alternatively, $\shY$ can be constructed using a graded ring of
theta functions, following \cite{GHKS}.

Here is a brief summary of the construction.
The central fibre $Y_0$ can be described as
\[
Y_0=\bigcup_{\sigma\in\P_{\max}}\PP_{\sigma}
\]
where $\P_{\max}$ denotes the maximal cells of $\P$ and $\PP_{\sigma}$
is the toric variety (projective if $\sigma$ is compact) determined by
the polyhedron $\sigma$. These toric varieties are glued together in a
manner reflecting the combinatorics of $\P$: if $\sigma_1\cap\sigma_2
=\tau$, then the strata $\PP_{\tau}\subseteq\PP_{\sigma_1}$, $\PP_{\tau}
\subseteq\PP_{\sigma_2}$ are identified.

Local models for the $k^{th}$ order deformation of $Y_0$ are 
determined by the function $\varphi$.  A key point of the
construction involves an invariant description for the local models,
which we explain here. The function $\varphi$, defined on an open
cover $\{U_i\}$ by single-valued functions $\varphi_i:U_i\rightarrow
Q^{\gp}_\RR$, determines an extension of $\Lambda$ by
$\underline{Q}^{\gp}$, the constant sheaf with coefficients in
$Q^{\gp}$. Indeed, on $U_i\cap B_0$, this extension will split as
$\Lambda|_{U_i}\oplus\underline{Q}^{\gp}$, and on the overlap, 
$(m,r)$ as a section of $\Lambda|_{U_i}\oplus\underline{Q}^{\gp}$
is  identified on $U_i\cap U_j$ with
$(m,r+d(\varphi_j-\varphi_i)(m))$ as a section of
$\Lambda|_{U_j}\oplus\underline{Q}^{\gp}$, interpreting
$d(\varphi_j-\varphi_i)\in \Hom(\Lambda|_{U_i}, Q^{\gp})$.  We then
have an exact sequence
\begin{equation}
\label{basicexactseq}
0\rightarrow \underline{Q}^{\gp}\rightarrow \shP\rightarrow \Lambda\rightarrow
0
\end{equation}
on $B_0$. We write the map $\shP\rightarrow\Lambda$ as $m\mapsto
\bar m$. After choosing a representative $\varphi_i$ of $\varphi$ in
a neighbourhood of a point $x\in B_0$, the stalk $\shP_x$ is
identified with  $\Lambda_x\oplus Q^{\gp}$.  There is a fan
$\Sigma_x=\{T_x\sigma\,|\,x\in\sigma\in \P\}$  (of not-necessarily
strictly convex cones), where $T_x\sigma$ denotes the tangent wedge
to $\sigma$ at $x$. This allows us to define  a convex PL function
$\varphi_x:|\Sigma_x| \rightarrow Q^{\gp}_\RR$ whose slope on
$T_x\sigma$ coincides with the slope of $\varphi_i|_{\sigma}$. We
then set
\begin{equation}
\label{Pvdef}
P_x:=\{(m,q)\,|\, m\in \Lambda_x\cap |\Sigma_x|, q\in Q^{\gp}, 
q-\varphi_x(m)\in Q\}
\subseteq\shP_x
\end{equation}
While this definition as described inside of $\Lambda_x\oplus Q$
depends on the choice of representative, in fact it is independent
of this choice when viewed as a submonoid of
$\shP_x$.

Note that $Q$ acts naturally on $P_x$, giving $\kk[P_x]$ a
$\kk[Q]$-algebra structure.  For a vertex $v$, we can now view
$\Spec \kk[P_v]/\fom^{k+1}$ as a local model for the $k^{th}$ order
deformation of $Y_0$ in a neighbourhood of the stratum of $Y_0$
corresponding to $v$. In addition, the local system $\shP$ gives a
method of defining parallel transport of monomials.

Let us describe certain aspects of this construction for our local
mirror symmetry example. Using the fan structure given by
\eqref{fanstructureatv},  we can describe the monoid $P_v\subseteq
\shP_v$ as $\{(m,r,s)\,|\, s- \varphi_v(m,r)\in Q\}\subseteq
M\oplus\ZZ\oplus Q^{\gp}$ using the identifications 
\begin{equation}
\label{canid}
\shP_v\cong\Lambda_v\oplus Q^{\gp}\cong 
M\oplus\ZZ\oplus Q^{\gp}
\end{equation} 
induced first by the representative $\varphi_v$ at $v$ and second
by the affine coordinate chart on $U_v$. In particular, for the 
purposes of the discussion below, we can describe the most relevant
part of $P_v$ as follows. First, we choose the representative 
$\varphi_v$ by choosing the linear function $\widetilde{\psi}_v$ 
to be $(0,\bar\psi(v),0)\in (N\oplus\ZZ^2)\otimes_{\ZZ} Q^{\gp}$,
with $\bar\psi=\psi|_{\sigma\times\{1\}}$.
Let $\bar P_v\subseteq P_v$ be the
submonoid consisting of $m\in P_v$ with $\bar m$ tangent to
$\sigma$. Then $\bar P_v$ is naturally described in terms of 
$\bar\psi$. Indeed, consider the
convex hull of the graph of $\bar\psi$,
\[
\Xi_{\bar\psi}:=
\{(m,0, s) \,|\, m\in\sigma, s - \bar\psi(m)\in Q\}
\subseteq M_{\RR}\oplus\RR \oplus Q^{\gp}_\RR,
\]
an unbounded polyhedron with vertices mapping to vertices of $\ol\P$
under the projection $M_\RR\oplus\RR\oplus Q_\RR^\gp\to M_\RR$. Then
we can identify $\bar P_v$ with the integral points in the tangent
wedge of $\Xi_{\bar\psi}$ at $(v,0,\bar\psi(v))$.

We also note that under the identification \eqref{canid} of 
$\shP_v$, the monodromy of $\Lambda$ described in
\eqref{monodromyformula} lifts to a monodromy transformation of
$\shP_v$ given by
\begin{align}
\label{monodromyformula2}
\begin{split}
T_{vv'}:\shP_v & \rightarrow \shP_v\\
(m,r,q) & \rightarrow (m+r(v-v'),r, q+ r(\bar\psi(v)-\bar\psi(v')))
\end{split}
\end{align}

The key additional (and usually most complex) ingredient for
constructing $\shY_k$ is a \emph{structure} $\scrS$. A structure
encodes data about how certain forms of these local models are glued
together. We will explain this structure in our example, but not go
into too much detail. A more detailed explanation for how this works
is given in the expository paper \cite{invitation}.

The structure takes a particularly simple form here. In general, a
structure is a collection of walls, polyhedral cells in $B$ of
codimension one each contained in a cell of $\P$ carrying the
additional data of certain formal power series. In
\cite{affinecomplex}  we distinguish a special sort of wall, namely
those contained in codimension one cells, and call them
\emph{slabs}. They tend to have a different behaviour. In the case
at hand, only slabs appear, and these cover $\sigma$.  The functions
attached to the slabs are determined from the monodromy around the
discriminant locus $\Delta$.

In this example, the slabs are the sets $\tau\times\{(1,0)\}$ for
$\tau\in\oP_{\max}$. For a slab $\fob$,  associated to any point
$x\in \fob\setminus \Delta$ is a formal power series
$f_{\fob,x}=\sum_{m\in P_x} c_mz^m$. This should only depend on the
connected component of $\fob\setminus \Delta$ containing $x$, so
there is in fact one such expression for each vertex $v$ of $\tau$,
and we can write $f_{\fob,v}=\sum_{m\in P_v} c_mz^m$. Furthermore,
$c_m\not=0$ implies $\bar m$ is tangent to $\fob$, so in fact the
sum is over $m\in \bar P_v$. 

The series $f_{\fob,v}$ are completely determined by a number of simple
properties. This follows in the case under consideration from having
chosen $\overline{\P}$ to consist of \emph{standard} simplices.
In what follows we will want to compare
$f_{\fob,v}$ with $f_{\fob,v'}$ for different vertices $v,v'$ of
$\fob$. To do so, we use parallel transport in $\shP$ from $v$ to $v'$.
Given the identification $\shP_v$ with $M\oplus\ZZ\oplus Q^{\gp}$ used
to give the formula \eqref{monodromyformula2} and noting that only
monomials of the form $z^{(m,0,p)}$ can appear in $f_{\fob,v}$, we see
that the particular path chosen between $v$ and $v'$ is irrelevant.

We can now state the conditions determining the $f_{\fob,v}$:
\begin{enumerate}
\item The constant term of each $f_{\fob,v}$ is $1$.
\item If $v$ and $v'$ are adjacent vertices of $\fob$, 
then the corresponding slab functions are related by
\begin{equation}
\label{slabrelations}
f_{\fob,v'}=z^{(v-v',0,\bar\psi(v)-\bar\psi(v'))}f_{\fob,v}.
\end{equation}
Here the equality makes sense after parallel transport of the exponents from
$v$ to $v'$ in the local system $\shP$, and $(v-v',0,\bar\psi(v)-
\bar\psi(v'))\in P_{v'}$ using the identification of $\shP_{v'}$ given
by \eqref{canid}.
\item $\log f_v$ contains no terms of the form $z^q$ for $q\in 
Q\setminus \{0\}$. Here we view $Q\subseteq P_v$ via the natural
inclusion $Q^{\gp}\subseteq \shP_v$.
\item If $v$ lies in slabs $\fob, \fob'$, then $f_{\fob,v}=f_{\fob',v}$.
\end{enumerate}

Item 1 is a normalization which originated in \cite{logmirror1}, Def.\ 4.23. 
However, we shall see its enumerative importance in \S \ref{Sect: Enumerative}.
Item 2 is the crucial point of
slabs: they allow us to define parallel transport of monomials
through slabs in a way which cancels the effects of monodromy. We
shall say more about this shortly. The condition 3 is interpreted by
writing $f_v=1+\cdots$ and using the Taylor expansion for
$\log(1+x)=\sum_{i=1}^{\infty} (-1)^{i+1}x^i/i$. This can be
interpreted inside some suitably completed ring. After expanding out
each expression $(\cdots)^i$, one demands that no monomials of the
form $z^q$ appear for any $q\in Q\setminus\{0\}$.  Finally, 4 tells
us how expressions propagate across $\sigma\times\{1\}$.

To see the significance of the second condition,  let $w\in
\Int(\tilde\sigma_1)$, $w'\in \Int(\tilde\sigma_2)$. Suppose we want
to compare monomials defined at $w$ (that is, monomials with
exponent in $P_w$) with monomials defined at $w'$ (that is,
monomials with exponent in $P_{w'}$). If we parallel transport from
$\shP_w$ to $\shP_{w'}$, the result depends on the path. For
example, let $v,v'$ be adjacent vertices of $\tau\in\oP_{\max}$. Let
$T_v$, $T_{v'}$ denote parallel transport in $\Lambda$ from $w$ to
$w'$ via the vertices $v$ and $v'$ respectively. Then from
\eqref{monodromyformula2}, it follows that for
$(m,r,q)\in\shP_w=M\oplus\ZZ\oplus Q$,
\[
T_{v'}(m,r,q)-T_v(m,r,q)=\big(r(v-v'),0,r(\bar\psi(v)-\bar\psi(v')\big).
\]
For convenience, we can identify $\shP_{w}$ and $\shP_{w'}$ with
$\shP_v$ so that $T_v$ is the identity. This difference between
$T_v$ and $T_{v'}$ creates problems for comparing the rings
$\kk[P_w]$ and $\kk[P_{w'}]$.  However, we can follow the rule that
if we wish to transport a monomial  $z^{(m,r,q)}$ along a path
between $w$ and $w'$ which crosses a slab $\fob$ in a  connected
component of $\fob\setminus\Delta$ containing a vertex $v$, we apply
an automorphism 
\begin{equation}
\label{transportthroughv}
z^{(m,r,q)} \mapsto z^{(m,r,q)}f_{\fob,v}^{-r}.
\end{equation}
Here $r$ represents the result of projecting
$\overline{(m,r,q)}=(m,r)$ via the projection
$\pi:\Lambda_v\rightarrow\ZZ$ obtained by dividing out by the
tangent space to the slab. If instead we pass through the slab
$\fob$ via the connected component of $\fob\setminus\Delta$
containing $v'$, we get
\[
z^{(m,r,q)}\longmapsto z^{(m+r(v-v'),r,q+r(\bar\psi(v)-\bar\psi(v'))} 
f_{v'}^{-r}= z^{(m,r,q)} f_v^{-r},
\]
coinciding with \eqref{transportthroughv}. Here we use the above expression
for $T_v-T_{v'}$ and \eqref{slabrelations}. Hence we see that the
ambiguity produced by monodromy is resolved by the slab functions.

\begin{examples}
\label{slabfunctionexamples}
In the following examples, we express the various functions
$f_{\fob,v}$ as formal power series with exponents appearing in
$\bar P_v$, using the representation of $\bar P_v$ as  the integral
points of the tangent wedge of $\Xi_{\bar\psi}$ at
$(v,0,\bar\psi(v))$.

(1) Take $\sigma$ to be the interval $[-1,1]$ as in Figure
\ref{localA1example}, with $\P$ as given there. The monoid of convex
piecewise linear functions on $\Sigma$ is generated by the function
which takes the values $0$, $0$ and $1$ respectively at $(-1,1)$,
$(0,1)$ and $(1,1)$. Thus we have $Q=\NN$, and the universal
piecewise linear function $\psi$ coincides with the above
generator. For a vertex $v$, with $\bar P_v\subseteq M\oplus 0
\oplus Q^{\gp}$, write $x=z^{(1,0,0)}$,  $t=z^{(0,0,1)}$, $t$ being
the generator of $\kk[Q]$. Then we have
\begin{align*}
f_{[-1,0],-1}= {} & 1+x+x^2t+xt,\\
f_{[-1,0],0}=f_{[0,1],0}= {} & 1+x^{-1}+xt+t,\\
f_{[0,1],1}= {} & 1+x^{-1}t^{-1}+x^{-2}t^{-1}+x^{-1}.
\end{align*}
Note that $\log f_{[-1,0],-1}$, $\log f_{[0,1],1}$ are clearly devoid
of pure powers of $t$ as any power, say, of  $x+x^2t+xt$ clearly
produces only terms with positive powers of $x$. On the other hand,
$f_{[-1,0],0}=(1+x^{-1})(1+xt)$, and taking logs we get $\log
(1+x^{-1})+\log (1+xt)$ which will again involve no pure $t$ power.
The $t$ term in $f_{[0,1],0}$ was  necessary to achieve this.

\medskip

(2) Take $\sigma$ to be as in Figure \ref{P2localfigure}. Again, 
the monoid of convex piecewise linear functions on the fan 
$\Sigma$ is generated by, say,
the function taking the values $0$ at $(0,0,1)$, $(1,0,1)$
and $(0,1,1)$ and the value $1$ at $(-1,-1,1)$. So again $Q=\NN$,
with the universal function $\psi$ agreeing with this generator.
Writing $x=z^{(1,0,0,0)}$,
$y=z^{(0,1,0,0)}$, $t=z^{(0,0,0,1)}$, 
it is easy to see that the terms of the slab function $f_{\fob,(0,0)}$
(independent of $\fob$ by the fourth condition) required by
conditions 1 and 2 are $1+x+y+tx^{-1}y^{-1}$. The normalization
condition forces us to add some additional terms:
\[
f_{\fob,(0,0)}=1+x+y+tx^{-1}y^{-1}+\sum_{k\ge 1} a_k t^k,
\]
where the $a_k$ are uniquely determined by the requirement that
\[
\sum_{i=1}^{\infty} (-1)^{i+1}
{(x+y+tx^{-1}y^{-1}+\sum_{k\ge 1} a_kt^k)^i\over i}
\]
contains no pure powers of $t$. This series in $t$ begins as
\[
-2t+5t^2-32t^3+286t^4-3038t^5+\cdots.
\]

(3) Let $\sigma$ be the convex hull of the points $(\pm 1,0)$, $(0,\pm 1)$ and
take $\P$ to be the star subdivision at the origin. Now the monoid $P$
of convex piecewise linear functions which are $0$ on $(0,0,1)$, $(1,0,1)$ and
$(0,1,1)$ is isomorphic to $\NN^2$, 
determined by the values $\alpha_1,\alpha_2$ of the function at generators
of the other two rays. Thus we can write $Q=\NN^2$,
$t_1=z^{(1,0)}\in \kk[Q]$,  $t_2=z^{(0,1)}\in\kk[Q]$. Using $x,y$ as
defined in the previous example, one can check that for any slab
$\fob$,
\[
f_{\fob,0}:= 1+x+y+t_1x^{-1}+t_2y^{-1}+t_1+t_2+3t_1t_2
+5t_1^2t_2+5t_1t_2^2
+\cdots.
\]
The additional terms represented by $\cdots$ give a power series in
$t_1,t_2$.
\end{examples}

We now describe the degeneration $\shY\rightarrow \Spec
\widehat{\kk[Q]}$ produced by the above data. In fact, it is not
difficult to do this in terms of equations, as follows. First,
define
\[
C(\Xi_{\bar\psi}):= \overline{\{((um,0,uq,u)\,|\,
(m,0,q)\in\Xi_{\bar\psi},\ 
u\in\RR_{\ge 0}\}}\subseteq M_{\RR}\oplus\RR\oplus Q^{\gp}_\RR\oplus\RR.
\]
Here the closure is necessary because $\Xi_{\bar\psi}$ is unbounded.
We then obtain a graded ring
\[
S_{\bar\psi}:=\kk[C(\Xi_{\bar\psi})\cap (M\oplus\ZZ\oplus Q^{\gp}\oplus\ZZ)]
\]
where the grading is given by the projection from  $M\oplus\ZZ\oplus
Q^{\gp}\oplus\ZZ$ onto the last copy of $\ZZ$.  Note the closure in
the definition of cone adds the cone $\{0\}\times\{0\}\times
\RR_{\ge 0}Q\times \{0\}$  to the set, so we see the degree $0$ part
of $S_{\bar\psi}$ is $\kk[Q]$. We can then complete, with
\[
\widehat S_{\bar\psi}:=S_{\bar\psi}\otimes_{\kk[Q]}\widehat{\kk[Q]}.
\]

It is then natural to think of the slab functions as being given by 
a single degree $1$ element of $\widehat S_{\bar\psi}$. Indeed,
given a vertex $v\in\oP$, we obtain from $f_{\fob,v}$ an element of
degree $1$ by multiplying all monomials of $f_{\fob,v}$ by
$z^{(v,0,\bar \psi(v),1)}$. It follows from \eqref{slabrelations}
that this is independent of the choice of $v$ and gives an element
$F\in\widehat{S}_{\bar\psi}$ of degree $1$. One can then show that
\begin{equation}\label{Eqn for shY}
\shY=\Proj\, \widehat S_{\bar\psi}[U,W]/(UW-z^q V_0 F).
\end{equation}
Here $U, W$ are of degree $1$, $V_0\in \widehat S_{\bar\psi}$ is the
element corresponding to $(0,0,0,1)$ (which lies in $\Xi_{\bar\psi}$
by the assumption that $0\in\sigma$ and $\psi$ has been chosen so that
$\bar\psi(0)=0$). The element $q\in Q$ is the element chosen in the definition
of $\widetilde\psi$ at the end of \S\ref{Sect: Construction}.
This can be shown in much the way the special case
discussed in \cite{invitation}, Example 5.2, being the case of
Examples \ref{slabfunctionexamples}, (2). Note that after localizing
at $z^q$, this family does not depend on the choice of $q$ up to isomorphism,
just as the choice of $q$ did not affect the polarization on the general 
fibres of $\shX\rightarrow\AA^1$.

The homogeneous coordinate ring of $\shY$  is generated in degree
$1$ by \emph{theta functions}, as explored in \cite{GHKS}. Each
point of $B(\ZZ)$ (the set of points of $B$ with integral
coordinates) corresponds to a generator of this ring as a
$\widehat{\kk[Q]}$-algebra. Explicitly, the integral points in this
example are the integral points of $\sigma$ and the apexes of the 
pyramids $\tilde\sigma_1$ and $\tilde\sigma_2$. If $v$ is an
integral point of  $\sigma$, then $z^{(v,0,\bar\psi(v),1)}\in
\widehat{S}_{\bar\psi}$ is the corresponding theta function. On the
other hand, the monomials $U$ and $W$ correspond to the two apexes.

This description of $\shY$ can be related to the more traditional
mirror to $X_{\Sigma}$ as described in \cite{CKYZ}. Here $\shY$ can
be decompactified by setting $V_0=1$, obtaining an open subset
$\shY^o$.  Passing to the generic fibre $\shY^o_{\eta}$ of
$\shY^o\rightarrow \Spec\widehat{\kk[Q]}$, we obtain a variety
defined over the field of fractions $K$ of $\widehat{\kk[Q]}$.  We
can describe $\shY^o$ as a subvariety of  $\AA^2\times
(N\otimes_{\ZZ}\Gm)$ over the field $K$ given by the equation
\begin{equation}
\label{tradmirror}
uw=z^qf_{\fob,0},
\end{equation}
where $\fob$ is any slab containing $0\in\sigma$. Here $u,w$ are
coordinates on $\AA^2$. Without the normalization condition, we
could take $f_{\fob,0}=\sum_{m\in\sigma\cap M}z^{(m,\bar\psi(m))}$,
which would lead to the mirror of $X_{\Sigma}$ being precisely that
given in \cite{CKYZ}. 

\begin{remark}
\label{Rem: Periods}
The crucial feature of the mirror family we have just described, as
opposed to the one given in \cite{CKYZ}, is that
the monomial coordinates on the base $\Spec \widehat{\kk[Q]}$ are
canonical in the sense of mirror symmetry. To describe this briefly,
we work over the field $\kk=\CC$, and assume 
that the power series
$f:=f_{\fob,0}$ is convergent in some analytic
neighbourhood $U$ of the zero-dimensional
stratum in $\Spec \CC[Q]$. Let $U^*=U \setminus \partial \Spec\CC[Q]$,
the complement of the union of toric divisors. Thus we can view $\shY^o$ as
giving an analytic family $\shY^o\rightarrow U^*$. 
We write $\shY^o_t$ for the fibre over $t\in U^*$. On such a fibre,
one has the holomorphic volume form on the fibres of $\cY^o\rightarrow
\Spec \widehat{\kk[Q]}$ given by
\[
\Omega=(2\pi i)^{-n-1}d\log u \wedge d\log x_1\wedge\cdots\wedge d\log x_n
\]
One then finds that there is a monodromy invariant cycle $\alpha_0\in
H_{n+1}(\cY^o_t,\ZZ)$ such that $\int_{\alpha_0}\Omega=1$, so that $\Omega$
is a normalized holomorphic form in the sense of mirror symmetry. Further,
if $q_1,\ldots,q_r\in Q^{\gp}$ are a basis for $Q^{\gp}$, one can find
(multi-valued) flat families of $(n+1)$-cycles $\alpha_1,\ldots,\alpha_r$ 
with $\int_{\alpha_i} \Omega=\log z^{q_i}$.
The key point of this calculation is to take the logarithmic derivative
of these period integrals and reduce the resulting integral to an integral
on the hypersurface $f_{\fob,v}=0$ in $N\otimes_{\ZZ}\Gm$. Via residues,
this is translated into an integral of the derivative of $\log f_{\fob,0}$
over various tori in $N\otimes_{\ZZ}\Gm$. The fact that these integrals
are then constant follows precisely from the normalization condition
on $f_{\fob,0}$.
\end{remark}

\section{Enumerative predictions}
\label{Sect: Enumerative}

So far we have seen two interpretations of the slab functions and
the normalization condition. The first came from the desire to write
down a correction to the patching of the naive toric models for the
mirror degeneration $\shY\to\Spec \widehat{\kk[Q]}$ in a way
consistent with local monodromy of the affine structure on $B$. We
discussed in \S\ref{Sect: Construction} how this condition along
with the normalization condition determines the slab functions
uniquely. Then in Remark \ref{Rem: Periods} we saw that the
normalization condition is responsible for making our families
canonically parametrized in the sense of mirror symmetry. Both of
these arguments concern the \emph{complex geometry} of the mirror
degeneration $\shY\to\Spec \widehat{\kk[Q]}$.

In the following two sections we will give two related
interpretations of normalized slab functions related to the
\emph{symplectic geometry} of the degeneration $\shX\to\AA^1$ of the
local Calabi-Yau variety $X_\Sigma$ we started with. The
interpretation supports the view that the degenerations constructed
by structures are indeed the ones expected from homological mirror
symmetry and open-closed string theory.

Since the completion of \cite{affinecomplex}, a clearer idea emerged as to
the precise meaning of structures. This picture has arisen from several
converging points of view:
(1) The heuristic correspondence between tropical Morse trees and Floer
homology emerging in discussions between us and Mohammed Abouzaid.
Some of these ideas were discussed in \cite{Clay} and \cite{thetasurvey}.
(2) Auroux's work \cite{Au} 
on $T$-duality on complements of anti-canonical divisors,
describing the complex structure on the SYZ dual of a Lagrangian fibration
using counts of Maslov index two disks. This has inspired
quite a bit of work, which is realising Fukaya's original
dream of correcting the complex structure of the mirror via counts of
holomorphic disks.
(3) \cite{GPS} made explicit the enumerative content of the key part
of the algorithm of \cite{KS} (or the two-dimensional version of
\cite{affinecomplex}). In particular, this established an enumerative
meaning for functions attached to walls of a structure.

Heuristically, one expects the following interpretation in the SYZ
picture of mirror symmetry. Suppose given a (special) Lagrangian
fibration $f:X\rightarrow B$ from a Calabi-Yau $X$, with the general fibre
being a torus. Consider Maslov index zero holomorphic disks
with boundary a fibre of $f$. For dimensional reasons the expectation
is that the set of points in $x\in B$ such that $f^{-1}(x)$ bounds a
Maslov index zero holomorphic disk is real codimension one in $B$,
forming a collection of walls. These walls should determine the structure
necessary to build the mirror to $X$, but one needs to attach functions
to the walls. Again, heuristically, these functions are expected to take
the shape, at a point $x\in B$ with $L=f^{-1}(x)$,
\begin{equation}
\label{nbetaformula}
\exp\left( \sum_{\beta\in \pi_2(X,L)\atop \partial\beta\not=0} k_{\beta}
n_{\beta}z^{\beta}
\right).
\end{equation}
Here the sum is over all relative homotopy classes $\beta$ such that
$\partial\beta\in\pi_1(L)$ is non-zero, 
$k_{\beta}$ is the index of $\partial\beta\in \pi_1(L)$ and
$n_{\beta}$ is some count of Maslov index zero disks with boundary
on $L$. This series should be defined as a
formal power series in some suitable ring. One can note that
as $x\in B$ varies, the groups $\pi_2(X,L)$ vary forming a local system
on $B_0$ (where $B_0=\{x\in B\,|\,\hbox{$f^{-1}(x)$ is non-singular}\}$).
This local system is analogous to the sheaf $\shP$ of \S 2, with the
exact sequence of homotopy groups
\[
\pi_2(L)=0 \rightarrow \pi_2(X)\rightarrow \pi_2(X,L)\rightarrow \pi_1(L)
\rightarrow \pi_1(X)
\]
being analogous to the exact sequence \eqref{basicexactseq}.

It is difficult to give exact definitions for the numbers $n_{\beta}$. 
There have been several approaches to dealing with this. For example,
Auroux \cite{Au} pioneered, in the case of an effective anti-canonical divisor,
the use of counts of Maslov index two disks to define holomorphic coordinates
which are then transformed by wall-crossing automorphisms as we cross 
walls in $B$ over which Maslov index zero disks live. 

A different approach,
using log geometry, originates in \cite{GPS}. There, working with Pandharipande,
we used relative Gromov-Witten invariants to make sense of the formula
\eqref{nbetaformula}. The situation there was effectively that of a 
rational surface with an anti-canonical divisor $D$, and the $n_{\beta}$ of
\eqref{nbetaformula} are replaced with counts of curves meeting the divisor
$D$ in one point. This was used for a general mirror symmetry construction
for such surfaces in \cite{GHK}.

It is interesting to note how these two points of view apply to the
case of local mirror symmetry considered in this paper. Auroux's point of
view was used effectively in a sequence of papers \cite{CLL}, \cite{CLT}, 
\cite{CCLT} 
to study the same local mirror symmetry situation as discussed in this paper.
Wall-crossing formulas for counts of Maslov index two disks 
are used to obtain what should be the same slab 
functions as discussed in this paper. The count of Maslov index two disks
is reduced to a closed Gromov-Witten invariant on a toric variety, which
can then be calculated via known mirror symmetry results. This
allows one to show
show that the slab functions defined using their counts give rise to
canonical coordinates just as our slab functions do.

On the other hand, generalising the idea of replacing holomorphic disks with
relative curves, one should be able to work with a certain kind of
logarithmic curve called a \emph{punctured curve}, the theory of which
is currently being developed in a joint project with Abramovich and Chen
\cite{ACGS}. These curves will live in the central fibre of the toric
degeneration $\cY\rightarrow \AA^1$ constructed in \S 2, and can be viewed
as a substitute for holomorphic curves with boundary in an algebro-geometric
context. Then \eqref{nbetaformula} can be used to define slab functions,
where now $n_{\beta}$ is a count of genus $0$ logarithmic curves with
one puncture.

We do not propose calculating the slab functions in this way. Rather,
we should be able to show that the slab functions defined in this way
satisfy the same determining properties that the slab functions of
\S 2 did. This is done by probing slabs by broken lines (see \cite{CPS},
\cite{GHK}, \cite{GHKS}) and interpreting
these enumeratively using a different type of punctured curve, roughly
corresponding to cylinders. These punctured curves play the same role that
Maslov index two disks play in the analysis of slab functions of
\cite{CLL}, \cite{CLT}, \cite{CCLT}. Crucially, we need to use the gluing
formula of \cite{ACGS} to relate broken lines and punctured curves.

\begin{figure}
\begin{center}
\input{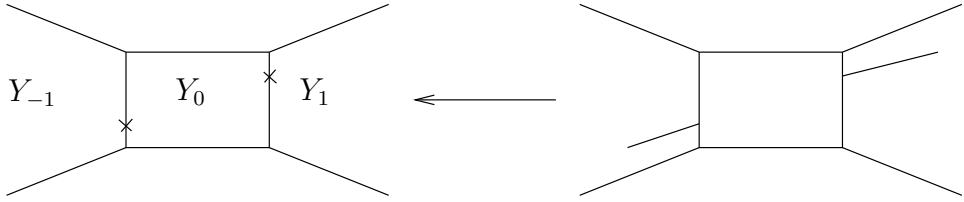}
\end{center}
\caption{The left-hand figure shows the five irreducible components of
$\shX_0$, with the three labelled components indexed by the vertices of
$\sigma$. Here $X_0\cong\PP^1\times\PP^1$ and $X_{-1}, X_1$ are isomorphic
to the blow-up of $\AA^2$ at a point.}
\label{A2exampleblowup}
\end{figure}

While the details of this approach will be given elsewhere, let us
demonstrate this using the simple example from Examples 
\ref{slabfunctionexamples}, (1). We depict
in Figure \ref{A2exampleblowup} 
the central fibre of the degeneration $\shX\rightarrow\AA^1$
constructed in \S 2 in this case. The total space $\shX$ has two ordinary
double points, situated on the singular locus of $\shX_0$, where the
map $\shX\rightarrow\AA^1$ is not normal crossings.
The inclusion $\shX_0\subseteq
\shX$ induces a log structure on $\shX_0$, but the log structure is
not well-behaved at the two points (not \emph{fine} in the sense of log
geometry). In particular, the theory of log Gromov-Witten invariants as
developed in \cite{JAMS}, \cite{AC}, \cite{Ch} 
cannot be used directly. While a theory of invariants
which can deal directly with this poorly behaved log structure is under
development, for the moment we will deal with it via a small resolution
of the ordinary double points. There are four choices of such resolutions,
one of which is shown on the right in Figure \ref{A2exampleblowup}.
These choices can be thought of in terms of the affine geometry
of the dual intersection complex $B$, 
with the resolutions corresponding to sliding the two
singularities of the affine structure along $\sigma$
to various choices of vertices.

\begin{figure}
\begin{center}
\input{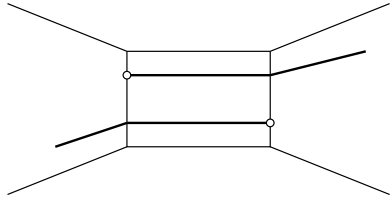}
\end{center}
\caption{The two punctured curves corresponding to holomorphic disks.
The curves include the exceptional divisors of the small resolution, and
the punctures are represented by the white circles.}
\label{A1curves}
\end{figure}

We have different slab functions $f_{\fob,v}$ for the vertices
$v=-1,0,1$. To identify the slab function at a vertex $v$ 
as a generating function, we choose a small resolution $\widetilde\shX
\rightarrow \shX$ so that the irreducible component indexed by $v$
remains toric. This effectively slides the singularities away from the
vertex. The resolution in Figure \ref{A2exampleblowup}
is used for the vertex $v=0$.
The slab function is given by \eqref{nbetaformula} where $n_{\beta}$
is a count of log curves of genus $0$ with one puncture mapping to the
boundary of the component $X_v$ indexed by $v$.
In Figure \ref{A1curves} we show the two obvious such curves for $v=0$. 
However, multiple covers of these
curves totally ramified at the puncture points are also possible,
and a $d$-fold cover will contribute with multiplicity $(-1)^{d+1}/d^2$.
The slab function is then, following \eqref{nbetaformula},
\[
(1+x^{-1})(1+xt)=\exp\left(\sum_{d=1}^{\infty} d \cdot {(-1)^{d+1}\over d^2}
x^{-d}+\sum_{d=1}^{\infty} d \cdot {(-1)^{d+1}\over d^2}(tx)^d\right),
\]
with the monomials $x^{-1}$ and $tx$ and their powers playing the role
of $z^{\beta}$.

\begin{figure}
\begin{center}
\input{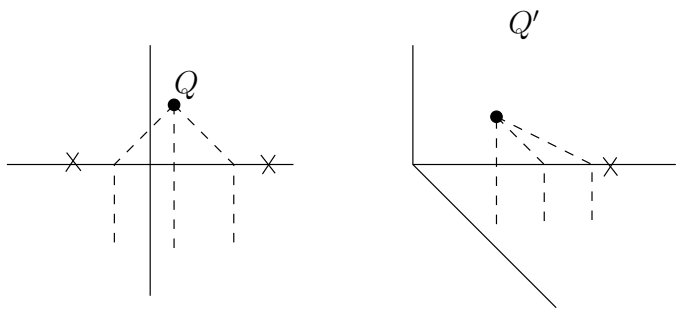}
\end{center}
\caption{There are four broken lines with endpoints $Q$ on the left,
two of which don't bend.
All have initial monomial $z^m$ with $\bar m=(0,-1)$. Once the broken
line crosses the slab, there are four possible attached monomials:
$z^m$, $tz^m$, $x^{-1}z^m$, $xtz^m$. The right-hand picture shows a
different choice of basepoint $Q'$, and there are again four broken lines.}
\label{brokenlines}
\end{figure}

To prove this formula without a direct calculation, we show the slab functions
defined by these counts satisfy conditions 1-4 of \S\ref{Sect: mirror}. 
Conditions~1 and
3 are obvious from \eqref{nbetaformula}, as the statement that only
monomials $z^{\beta}$ with $\partial\beta\not=0$ appear is analogous
to the statement that no terms of the form $z^q$ for $q\in Q\setminus
\{0\}$ appear inside the exponential. Condition 4 is automatic because
in this situation the slab function only depends on the vertex. It remains
to show condition 2, and we use broken lines for this, which can be reviewed
in \cite{thetasurvey}. A broken line is a piecewise linear path with monomials
$c_Lz^{m_L}$ attached to each domain of linearity, and the derivative of the
line in the domain $L$ is $-\bar m_L$. When the broken line crosses a wall,
we may change the monomial by applying the wall-crossing automorphism
\eqref{transportthroughv}
to the monomial and choosing a new monomial being one of the terms in the
expression obtained after applying this automorphism. In Figure 
\ref{brokenlines}, we consider germs of broken lines which
come vertically from below with initial monomial $z^m$ with $\bar m=(0,-1)
\in M_{\RR}\oplus\RR$. We define $\Lift_Q(m)$ to be the sum over all broken
lines ending at a basepoint $Q$ of the final attached monomials. 
Note that if $Q$ is near a vertex $v$
of $\sigma$, then in fact $\Lift_Q(m)=z^m f_{\fob,v}$. It is then not difficult
to show that \eqref{slabrelations} holds for all pairs of adjacent vertices
if and only if $\Lift_Q(m)$ is independent of $Q$
chosen above the slabs as in Figure \ref{brokenlines}.

Broken lines can be viewed as a purely combinatorial (tropical) way to
count holomorphic cylinders. But we can actually count 
logarithmic curves of genus $0$ with two punctures to emulate cylinders,
and there will be a correspondence between such twice-punctured logarithmic
curves and broken lines. Varying the basepoint can be achieved by varying
a point constraint for one of the punctures. 
The key point is that various ways of degenerating the point constraint
can lead to different broken lines with different endpoints. However, the
count of these punctured curves will be independent of the constraint.

To see this explicitly, let's look at the example of the straight line
in the left-hand diagram in Figure \ref{brokenlines} with attached
monomial $z^m$. To understand what happens when we move this broken line
through the singularity, it is helpful to move the singularity to the
vertex $0$ by using the small resolution depicted in Figure \ref{brokenlines2}.
Consider
the family of twice-punctured curves given by the vertical fibres of
$X_0$, the blowup of $\PP^1\times\PP^1$. 
Any curve in this family has a tropicalization
(see \cite{JAMS}, \S 3). The tropicalization of a general curve in this
family is just the vertical line through the singularity on the right-hand
side of Figure \ref{brokenlines2}. Combinatorially, this just indicates that
the curve intersects the upper and lower boundary divisors of $X_0$.
However, the family has two special members which are 
degenerate with respect to the log structure on the central fibre.
The tropicalization of the punctured curve when it falls into $X_0\cap X_1$
is the straight broken line depicted in Figure \ref{brokenlines2} 
to the right of the vertex.
If on the other hand we move the punctured curve to the left,
it becomes reducible, the union of $X_{-1}\cap X_0$ and the exceptional
curve of the small resolution. This curve tropicalizes to the tropical
curve depicted on the left, now with two vertices corresponding to the
two components. The bend is a consequence of gluing the once-punctured
curve with support the exceptional curve to the thrice-punctured curve
with support $X_0\cap X_{-1}$. The broken line is then a subset
of this tropical curve.

The point is that the two broken lines make the same contribution
to $\Lift_Q(m)$ as $Q$ varies because they can both be viewed
as counting the number of curves in the one-parameter family described
passing through some point in $X_0$. The point can degenerate into 
$X_{-1}\cap X_0$ or $X_0\cap X_1$, giving the two types of broken line
behaviour. Thus the invariance of $\Lift_Q(m)$ can be viewed as the
fact that these lifts are generating functions for counts of certain types
of punctured curves.

The key point for the argument is then to prove that broken lines really
calculate Gromov-Witten invariants of punctured curves. This shall be 
shown using a general gluing formula \cite{ACGS}. 

\begin{figure}
\begin{center}
\input{brokenlines2.pstex_t}
\end{center}
\caption{}
\label{brokenlines2}
\end{figure}

Note so far we have not actually calculated the Gromov-Witten invariants
of $X_{\Sigma}$. These should be extracted in the $B$-model from some
additional period integrals past the ones discussed in Remark
\ref{Rem: Periods}. A
significant challenge remaining is to give a tropical description for these
period integrals and Gromov-Witten invariants.

\section{Tropical disks and slab functions}

The picture of counting holomorphic disks and cylinders from
\S\ref{Sect: Enumerative} suggests an interpretation of the slab
functions in terms of tropical curves. In this section we give a
surprisingly simple interpretation of this sort. The arguments are
by algebraic manipulations of the slab functions. We are thus lead
to the challenge of interpreting the tropical counts in terms of the
counting of holomorphic disks on $X_\Sigma$.

We study the collection of slab functions at a vertex $v\in\ol\P$
with $v\in\Int\sigma$. By Condition~(4) of slab functions all the
$f_{\fob,v}$ for slabs $\fob$ containing $v$ agree. Dehomogenizing
\eqref{Eqn for shY} at $v$ we are thus left with the local model
$uw-ft=0$ for the mirror degeneration for some
$f\in\widehat{\kk[P]}$. Here $P=\bar P_v$ is a toric submonoid of $M\oplus
Q^\gp$ with $P^\times=\{0\}$ and the completion is with respect to
$P\setminus\{0\}$. Recall also the projection
\[
M\oplus Q^\gp\lra M,\quad m\longmapsto \bar m. 
\]

For example, for the mirror of local $\PP^2$ we had $Q=\NN$,
$P\subset\NN^3$ generated by $(1,0,0)$, $(0,1,0)$, $(-1,-1,1)$,
hence $\widehat{\kk[P]}= \kk[x,y,z] \lfor t\rfor/(xyz-t)$, and
\[
f=1+x+y+z-2t+5t^2-32t^3+286t^4-3038t^5+\cdots.
\]

In general we assume $f=1+\sum_{i=1}^r z^{m_i}+g$ with $\bar m_i\neq
0$ for all $i$ and $g=\sum_q \mbox{$b_q\cdot z^q$}\in
\widehat{\kk[Q]}$ taking care of the normalization
condition.\footnote{Note that by the universal nature of $Q$ the sum
over $z^{m_i}$ implicitly comprises a universal choice of
coefficients.}  Under this assumption we are going to give an
infinite product expansion
\[
f=\prod_{\{m\,|\, \bar m\neq 0\}} (1+a_m z^m)
\]
in $\widehat{\kk[P]}$, with each $a_m$ having an interpretation in
terms of tropical disks in $M_\RR$ with root weight $m$. Moreover,
each coefficient $b_q$ of $g$ has an interpretation in terms of
pointed tropical curves of genus zero.\footnote{If $\sigma$ has
several interior integral points the change of vertex
formula~\eqref{monodromyformula2} provides a non-trivial identity
between expressions labelled by different sets of tropical trees. It
would be interesting to give an interpretation of this formula
within the following discussion.}

To this end consider the following definition of the \emph{type of a
tropical disk}. A \emph{rooted tree} is a partially ordered finite
set with a unique maximal element, called the \emph{root vertex},
which is connected and without cycles when viewed as a graph. The
\emph{predecessors} of a vertex $v$ are the adjacent vertices that
are smaller than $v$. The minimal elements of a tree are called its
\emph{leaves}, so these are the elements without predecessors. We
require that there are no elements with only one predecessor. In
graph theory language this means that the interior vertices are at
least trivalent and the leaves are the unique univalent vertices.

We now define types of tropical trees weighted by elements of $P$.
Note that $Q$ can be identified with the submonoid $\{m\in P\,|\,
\bar m=0\}$ of $P$.

\begin{definition}\label{Def: tropical disks}
The \emph{type of a $P$-labelled tropical disk} is a rooted tree
$\Gamma$ with sets $V_\Gamma$ of vertices and $E_\Gamma$ of edges
along with a \emph{vertex-labeling} map
\[
w:V_\Gamma\lra P\setminus Q,\quad v\longmapsto m_v
\]
fulfilling the following conditions:
\begin{enumerate}
\item
For any non-leaf vertex $v\in V_\Gamma$ with predecessors
$v_1,\ldots,v_\ell$ the balancing condition
\[
m_v= m_{v_1}+\cdots+m_{v_\ell}
\]
holds.
\item
For any vertex $v$ the weights $m_1,\ldots,m_\ell$
of the adjacent predecessor vertices are pairwise distinct.
\end{enumerate}
By abuse of notation we suppress the labelling function in the notation
and write just $\Gamma$ for the type of a tropical disk. The set of
non-leaf vertices is denoted $\hat V_\Gamma$.

If we take the weight $m_\Gamma$ of the root vertex in
$Q$ rather than in $P\setminus Q$ and otherwise leave the definition
unchanged we arrive at the notion of \emph{type of $P$-labelled
pointed rational tropical curve}.
\end{definition}

Each type of tropical disk or rational tropical curve determines an
isotopy class of traditional tropical curves in $M_\RR$ with edges
labelled by lifts of the direction vector (an element of $M$) to
$P$, the labelling of the predecessor vertex. In the disk case one
may add another edge to force the balancing condition at the root
vertex.

The balancing condition for a tropical disk implies that the
labelling function is uniquely determined by its values on the leaf
vertices $v_1,\ldots,v_\ell$. In particular, for the weight of the
root vertex we have
\[
m_\Gamma=m_{v_1}+\cdots+ m_{v_\ell}.
\]

Let now $S=\{m_1,\ldots,m_r\}$ be the set of exponents $m$ occuring
in $f$ with $\bar m\neq0$. For $m\in P$ with $\bar m\neq 0$ denote by
$\scrT_m(S)$ the set of types of $P$-labelled tropical disks
$\Gamma$ with $m_\Gamma= m$ and with leaf labels in $S$.
Similarly $\scrR_q(S)$ denotes the set of types of $P$-labelled
pointed rational tropical curves with leaf labels in $S$ and $m_{\Gamma}=q$.

\begin{proposition}
For $f=1+\sum_{m\in S} z^m+g$ with $g=\sum_{q\neq 0} b_q z^q$
as above it holds
\begin{equation}
\label{Eqn: Tropical product expansion}
f=\prod_{\{m\,|\, \bar m\neq0\}} (1+a_m z^m)
\end{equation}
with
\[
a_m= \sum_{\Gamma\in\scrT_m(S)} (-1)^{|\hat V_\Gamma|}\quad
\text{and}\quad b_q= \sum_{\tilde\Gamma\in \scrR_q(S)}
(-1)^{|V_{\tilde\Gamma}|-1}.
\]
\end{proposition}

\proof
Expanding the infinite product in the statement and gathering
according to monomials yields
\begin{equation}\label{Eqn: Expansion}
\prod_{\{m\,|\, \bar m\neq0\}} (1+a_m z^m)=
\sum_{m\in P} \bigg( \sum_{\ell=1}^{\infty} 
\sum_{m=m_1+\cdots+m_\ell \atop \Gamma_i\in\scrT_{m_i}(S)}
(-1)^{|\hat V_{\Gamma_1}|}\cdots (-1)^{|\hat V_{\Gamma_\ell}|} \bigg) z^m.
\end{equation}
In this expansion $\ell$ is the number of $a_m$-terms in the
infinite product to be multiplied. Thus the third sum on the
right-hand side is over all decompositions $m=m_1+\cdots+m_\ell$ of
$m$ into $\ell$ \emph{pairwise distinct} summands in $P$. Recall
that $\hat V_\Gamma$ is the set of non-leaf vertices.
Fix $m$ with $\bar m\neq 0$ now and consider the coefficient of
$z^m$. Then for $\ell\ge 2$ any collection $\Gamma_1,\ldots,
\Gamma_\ell$ of types of tropical disks with $m=m_{\Gamma_1}+\cdots
+m_{\Gamma_\ell}$ can be merged into a new type of tropical disk
$\Gamma\in \scrT_m(S)$ by connecting the root vertex of each
$\Gamma_i$ by one edge to the root vertex $v_0\in V_\Gamma$. Thus the
root vertex of $\Gamma$ is $\ell$-valent with adjacent predecessor
trees $\Gamma_1,\ldots,\Gamma_\ell$. Now this merged tree $\Gamma$
contributes to the coefficient of $z^m$ as one term for $\ell=1$.
Since the vertices of $\Gamma$ other than the root vertex are in
bijection with the vertices of $\Gamma_1,\ldots,\Gamma_\ell$ it
holds
\[
(-1)^{|\hat V_\Gamma|}=
-(-1)^{|\hat V_{\Gamma_1}|}\cdots (-1)^{|\hat V_{\Gamma_\ell}|}.
\]
Thus each term with $\ell\ge 2$ in the sum of the right-hand side of
\eqref{Eqn: Expansion} cancels with one term for $\ell=1$.
Conversely, if the root vertex of the type of a tropical disk
$\Gamma$ has valency $\ell\ge 2$ then $\Gamma$ is obtained by this
merging procedure. On the right-hand side of \eqref{Eqn: Expansion}
we are thus left only with those $m$ with $\bar m=0$ and in addition
with those trees with only one vertex. The latter condition means
that the root vertex is also a leaf vertex. These terms yield the
sum $\sum_{m\in S} z^m$. The terms with $\bar m=0$ define a power
series $1+h\in\widehat{\kk[Q]}$. We have thus shown
\[
\prod_{\{m\,|\, \bar m\neq0\}} (1+a_m z^m)= 1+\sum_{i=1}^r z^{m_i}+h,
\]
with $h\in \widehat{\kk[Q]}$. Since the left-hand side of this
equation is clearly normalized we see that $h=g$. Tropically, the
coefficient of $z^q$ in $g$ is the weighted sum of types of
$P$-labelled pointed rational tropical curves, with the marked point
(of valency $\ell\ge 2$) the merging point of $\ell$ tropical trees
$\Gamma_1,\ldots,\Gamma_\ell$. The balancing condition of an
underlying tropical curve at the marked point is the statement
\[
\bar m_{\Gamma_1}+\cdots+\bar m_{\Gamma_\ell}
= \bar m=0.
\]
\qed
\begin{figure}
\begin{center}
\input{trees.pstex_t}
\end{center}
\caption{}
\label{trees}
\end{figure}

For $f=1+x+y+z+g$ the expansion up to order $4$ is
\begin{eqnarray*}
&&(1+x)(1+y)(1+z)(1-xy)(1-yz)(1-xz)(1+x^2y)\\
&&\cdot(1+xy^2)\ldots(1+yz^2)(1-x^2y^2)(1-y^2z^2)(1-x^2z^2)\\
&&\cdot(1-x^2yz)(1-xy^2z)(1-xyz^2)(1-xz^3)\ldots(1-yz^3)
\end{eqnarray*}
Figure~\ref{trees} shows the tropical trees contributing to the coefficient
$-1=(-1)^3+(-1)^3+(-1)^2$ of $x^2y^2$. Note that many labelled trees
with four leaves are ruled out because of the third condition in
Definition~\ref{Def: tropical disks} that no two predecessor
subtrees at some vertex be isomorphic.

We finish this section with two remarks on a possible enumerative
interpretation of the expansion in terms of tropical disks and trees.
First, according to \eqref{nbetaformula} we should write the product
expansion \eqref{Eqn: Tropical product expansion} in exponential
form. Indeed, we can also write
\[
f=\exp\left(\sum_{\{m\,|\,\bar m\neq 0\}}
\sum_{\Gamma\in\tilde\scrT_m(S)} \frac{(-1)^{|\hat
V_\Gamma|}}{|\Aut(\Gamma)|}z^m\right).
\]
Here the sum is over the space $\tilde\scrT_m(S)$ of tropical disks
with the stability condition Definition~\ref{Def: tropical disks},2
dropped. Expanding $\exp$ in a Taylor series the proof is largely
the same as the one given, with extra care taken concerning automorphisms.

Second, in log Gromov-Witten theory the log structure on the moduli
space only depends on the type of tropical curve associated to a
stable log map \cite{JAMS}. It is tempting to believe in a
formulation of the counting problem by a symmetric obstruction
theory \cite{BehrendFantechi} on a moduli space with a log structure
stratified by types of tropical disks, with each stratum
contributing $(-1)^{|\hat V_\Gamma|}/|\Aut(\Gamma)|$.


\end{document}